\documentclass{ifacconf}

\usepackage{graphicx}      
\usepackage{natbib}        
\usepackage{amsmath}
\usepackage{amssymb}

\newtheorem{theorem}{\indent Theorem}
\newtheorem{remark}{\indent Remark}

\begin{document}

\def \eps{\varepsilon}

\begin{frontmatter}

\title{Batch Data Processing \\and Gaussian Two-Armed Bandit
\thanksref{footnoteinfo}}

\thanks[footnoteinfo]{This work was supported in part by the Project Part of the State
Assignment in the Field of Scientific Activity by the Ministry of
Education and Science of the Russian Federation, project no.
1.949.2014/K.}

\author[First]{Alexander V. Kolnogorov}

\address[First]{Yaroslav-the-Wise Novgorod State University,
   Velikiy Novgorod, 173003 Russia, (e-mail: Alexander.Kolnogorov@novsu.ru).}

\begin{abstract}                
We consider the two-armed bandit problem as applied to data
processing if there are two alternative processing methods
available with different a priori unknown efficiencies. One should
determine the most effective method and provide its predominant
application. Gaussian two-armed bandit describes the batch, and
possibly parallel, processing  when the same methods are applied
to sufficiently large packets of data and accumulated incomes are
used for the control. If the number of packets is large enough
then such control does not deteriorate the control performance,
i.e. does not increase the minimax risk. For example, in case of
50 packets the minimax risk is about 2\% larger than that one
corresponding to one-by-one optimal processing. However, this is
completely true only for methods with close efficiencies because
otherwise there may be significant expected losses at the initial
stage of control when both actions are applied turn-by-turn. To
avoid significant losses at the initial stage of control one
should take initial packets of data having smaller sizes.
\end{abstract}

\begin{keyword}
two-armed bandit problem, stochastic robust control, minimax and
bayesian approaches, batch processing.
\end{keyword}

\end{frontmatter}

\section{Introduction}

We consider the following setting of the two-armed bandit problem
(see, e.g. \cite{BF}, \cite{PS}) which is also well-known as the
problem of expedient behavior in a random environment (see, e.g.
\cite{Tsetlin}, \cite{Varsh}) and the problem of adaptive control
in a random environment (see, e.g. \cite{Sragovich}, \cite{NP}).
Let $\xi_n$, $n=1,\dots,N$ be a controlled random process which
values are interpreted as incomes, depend only on currently chosen
actions $y_n$ ($y_n \in \{1,2\}$) and are normally distributed
with probability densities $f(x|m_\ell)$ if $y_n=\ell$,
$\ell=1,2$, where
\begin{gather*}
f(x|m)=(2\pi)^{-1/2} \exp \left\{-(x-m)^2/2\right\}.
\end{gather*}
It is the so-called Gaussian (or Normal) two-armed bandit. It can
be completely described by a vector parameter $\theta=(m_1, m_2)$.
The goal is to maximize the total expected income. To thus end,
one should determine the action corresponding to the largest value
of $m_1$, $m_2$ and provide its predominant application.
\par
Let's explain why Gaussian two-armed bandit is considered. We
investigate the problem as applied to control of data processing
if there are two alternative processing methods available with
different a priori unknown efficiencies. Let $T=NM$ items of data
be given which may be processed by either of the two alternative
methods. Processing may be successful ($\zeta_t = 1$) or
unsuccessful ($\zeta_t = 0$). The goal is to maximize  the total
expected number of successfully processed items of data.
Probabilities of successful and unsuccessful processing depend
only on chosen methods (actions), i.e. $\Pr(\zeta_t =
1|y_t=\ell)=p_\ell$, $\Pr(\zeta_t = 0|y_t=\ell)=q_\ell$,
$\ell=1,2$. Assume that $p_1$, $p_2$ are close to $p$ ($0<p<1$).
We partition all data items into $N$ packets each containing $M$
data items. For data processing in each packet we use the same
method. Note that data in the same packet may be processed in
parallel. For control we use the values of the process
$\xi_n=(DM)^{-1/2}\displaystyle{\sum_{t=(n-1)M+1}^{nM}} \zeta_t$,
$n=1,\dots,N$ with $D=p(1-p)$. According to the central limit
theorem distributions of $\xi_n$, $n=1,\dots,N$ are close to
Gaussian and their variances are close to unity just like in
considered setup.
\par
\begin{remark}
Parallel control in the two-armed bandit problem was first
proposed for treating a large group of patients by either of the
two alternative drugs with different unknown efficiencies.
Clearly, the doctor cannot treat the patients sequentially
one-by-one. Say, if the result of the treatment will be manifest
in a week and there is a thousand of patients, then one-by-one
treatment would take about twenty years. Therefore, it was
proposed to give both drugs to sufficiently large groups of
patients and then the most effective one to give to the rest of
them. As the result, the entire treatment takes two weeks. The
discussion and bibliography of the problem as applied to medical
trials can be found, for example, in~\cite{Lai, Cheng}.
\end{remark}
\par
Control strategy $\sigma$ at the point of time $n$ assigns a
random choice of the action $y_n$ depending on the current history
of the process, i.e. responses $x^{n-1}=x_1,\dots,x_{n-1}$ to
applied actions $y^{n-1}=y_1,\dots,y_{n-1}$:
\begin{gather*}
\sigma_{\ell}(y^{n-1},x^{n-1})=\Pr(y_n=\ell|y^{n-1},x^{n-1}),
\end{gather*}
$\ell=1,2$. The set of strategies is denoted by $\Sigma$.
\par
Recall that the goal is to maximize (in some sense) the total
expected income. Therefore, if parameter $\theta$ is known then
the optimal strategy should always apply the action corresponding
to the largest value of $m_1$, $m_2$. The total expected income
would thus be equal to $N (m_1 \vee m_2)$. If the parameter is
unknown then the loss function
\begin{equation}\label{c1}
L_N(\sigma,\theta)=N (m_1 \vee m_2)-
\mathbb{E}_{\sigma,\theta}\left(\sum_{n=1}^N  \xi_n\right)
\end{equation}
is equal to expected losses of total income with respect to its
maximal possible value. Here $\mathbb{E}_{\sigma,\theta}$ denotes
the mathematical expectation calculated with respect to measure
generated by strategy $\sigma$ and parameter $\theta$. The set of
parameters is assumed to be the following
\begin{gather*}
\Theta=\{\theta:|m_1-m_2| \leq 2C\},
\end{gather*}
where $0<C< \infty$. Here restriction $C < \infty$ ensures the
boundedness of the loss function on $\Theta$.
\par
According to the minimax approach the maximal value of the loss
function on the set of parameters $\Theta$ should be minimized on
the set of strategies $\Sigma$. The value
\begin{equation}\label{c2}
R^M_N(\Theta)=\inf_{\Sigma} \sup_{\Theta} L_N(\sigma,\theta)
\end{equation}
is called the minimax risk and corresponding strategy $\sigma^M$
is called the minimax strategy. Note that if strategy $\sigma^M$
is applied then the following inequality holds
\begin{equation*}
L_N(\sigma^M,\theta) \le R^M_N(\Theta)
\end{equation*}
for all $\theta \in \Theta$ and this implies robustness of the
control.
\par
The minimax approach to the problem was proposed
in~\cite{Robbins} and caused a considerable interest to it. The
classic object of most of arisen articles was the so-called
Bernoulli two-armed bandit which can be described by distribution
\begin{gather*}
\Pr(\xi_n =1|y_n=\ell)=p_\ell,\quad \Pr(\xi_n =0|y_n=\ell)=q_\ell,
\end{gather*}
$p_\ell+q_\ell=1$, $\ell=1,2$. It can be described by a parameter
$\theta=(p_1,p_2)$ with the set of values $\Theta=\{\theta:0\leq
p_\ell\leq 1; \ell=1,2\}$. It was shown in \cite{FZ} that explicit
determination of the minimax strategy and minimax risk is
virtually impossible already for $N\ge 5$. However, an asymptotic
minimax theorem was proved in \cite{Vogel} using some indirect
techniques. This theorem states that the following estimates hold
as $N \to \infty$:
\begin{equation}
0.612  \leq (DN)^{-1/2}R^M_N(\Theta) \leq 0.752, \label{c3}
\end{equation}
where $D=0.25$ is the maximal variance of one-step income.
Presented here the lower estimate was obtained in \cite{Bather}.
The maximal value of expected losses corresponds to
$|p_1-p_2|\approx 3.78 (D/N)^{1/2}$ with additional requirement
that $p_1$, $p_2$ are close to $0.5$.
\begin{remark}
There are some different approaches to robust control in the
two-armed bandit problem, see, e.g. \cite{NP, Lugosi, UN, GNS}. In
these articles, another ideas like stochastic approximation method
and mirror descent algorithm are used for the control. The order
of the minimax risk for these algorithms is $N^{1/2}$ or close to
$N^{1/2}$.
\end{remark}
\par
Another very popular approach to the problem is a Bayesian one.
Let $\lambda(\theta)=\lambda(m_1,m_2)$ be some prior probability
density. The value
\begin{equation}\label{c4}
R^B_N(\lambda)=\inf_{\Sigma} \int_{\Theta}
L_N(\sigma,\theta)\lambda(\theta) d\theta
\end{equation}
is called the Bayesian risk and corresponding optimal strategy is
called the Bayesian strategy. Bayesian approach allows to find
Bayesian strategy and risk by solving a recursive Bellman-type
equation. Minimax risk~\eqref{c2} and Bayesian risk~\eqref{c4} are
related by the main theorem of the theory of games as follows:
\begin{equation}\label{c5}
R^M_N(\Theta)=R^B_N(\lambda_0)= \sup_{\lambda} R^B_N(\lambda),
\end{equation}
where $\lambda_0$ is called the worst-case prior distribution.
\par
The goal of this paper is to present the approach based on the
main theorem of the theory of games. This approach allows to
determine minimax strategy and minimax risk explicitly by solving
appropriate Bellman-type recursive equation and finding the
worst-case prior distribution according to~\eqref{c5}. This allows
to evaluate the control performance. In particular, it turned out
that in case of close mathematical expectations $m_1$, $m_2$
batching of data almost does not enlarge the maximal expected
losses if the number of packets is large enough, e.g. if the
number of packets is 50 or larger. Therefore, say $50000$ items of
data may be processed in 50 steps by packets of 1000 data with
almost the same maximal losses as if the data were processed
optimally one-by-one. To be more precise, the maximal expected
losses in case of batch processing in 50 steps are about 2\%
larger than in case of optimal one-by-one processing. However, in
case of distant expectations there may be large expected losses at
the initial stage of control when actions are applied
turn-by-turn. To reduce the losses at the initial stage, one
should reduce corresponding sizes of packets. The example is given
in Section~\ref{CIS}
\par
The structure of the paper is the following. In section~\ref{MT}
we present the Bellman-type recursive equation which allows to
determine explicitly Bayesian strategy and risk for any prior
distribution. In section~\ref{WCP} properties of the worst-case
prior distribution are investigated and this allows to simplify
the recursive Bellman-type equation significantly. In
section~\ref{IRE} we obtain invariant recursive Bellman-type
equation with unity control horizon and  its limiting description
by the second order partial differential equation. In
section~\ref{CIS} we find minimax risks numerically.
Section~\ref{Con} contains a conclusion. Note that some results
are presented in \cite{Koln1}, \cite{Koln2}, \cite{Koln3}. Here we
combine and compare them.

\section{Recursive equation for
determination of Bayesian strategy and risk} \label {MT}

Bayesian strategy and risk can be calculated recursively. Let
history of control up to the point of time $n$ be described by
$(X_1,n_1,X_2,n_2)$. Here $n_1$, $n_2$ are total numbers of
applications of both actions ($n_1+n_2=n$) and $X_1$, $X_2$ are
corresponding total incomes. Let $X_\ell=0$ if $n_\ell=0$. Denote
by $f_D(x|m)=(2\pi D)^{-1/2} \exp \left\{-(x-m)^2/(2D)\right\}$
the Gaussian probability density. The posterior distribution
density is thus equal to
\begin{equation*}
\begin{array}{ll}
\lambda(m_1,m_2|X_1,n_1,X_2,n_2)\\
=\displaystyle{\frac{f_{n_1}(X_1|n_1 m_1)f_{n_2}(X_2|n_2 m_2)
\lambda(m_1,m_2)} {p(X_1,n_1,X_2,n_2)}}
\end{array}
\end{equation*}
with
\begin{gather}\label{c6}
p(X_1,n_1,X_2,n_2)\\=
\displaystyle{\iint\limits_{\Theta}}f_{n_1}(X_1|n_1
m_1)f_{n_2}(X_2|n_2 m_2) \lambda(m_1,m_2)dm_1dm_2 \nonumber
\end{gather}
If it is assumed additionally that $f_n(X|nm)=1$ at $n=0$ then
this expression holds true if $n_1=0$ and/or $n_2=0$ as well.
\par
In the sequel, we consider strategies which apply each chosen
action $M$ times. For the sake of simplicity we assume that $N$ is
a multiple of $M$. If incomes arise sequentially one-by-one, these
strategies allow to switch actions more rarely. If incomes arise
by packets, these strategies allow their parallel processing.
Denote by $R^B_{N-n}(\lambda;X_1,n_1,X_2,n_2)$ Bayesian risk at
the latter $(N-n)$ steps calculated with respect to the posterior
distribution density $\lambda(m_1,m_2|X_1,n_1,X_2,n_2)$. Let
$x^+=\max(x,0)$. Then
\begin{equation}
R^B_{N-n}(\cdot)=\min(R^{(1)}_{N-n}(\cdot),R^{(2)}_{N-n}(\cdot)),
\label{c7}
\end{equation}
where $R^{(1)}_{0}(\cdot)=R^{(2)}_{0}(\cdot)=0$ if $n_1+n_2=N$,
\begin{equation}
\begin{array}{lll}
R^{(1)}_{N-n}(\lambda;X_1,n_1,X_2,n_2)=
\displaystyle{\iint\limits_{\Theta}} \left(\right. M (m_2-m_1)^+
\\
\quad + \mathbb{E}_{x}^{(1)}
R^B_{N-(n+M)}(\lambda;X_1+x,n_1+M,X_2,n_2) \left.\right) \\
\qquad  \times  \lambda(m_1,m_2|X_1,n_1,X_2,n_2) dm_1 dm_2,
\end{array} \label{c8}
\end{equation}
\begin{equation}
\label{c9}
\begin{array}{lll}
R^{(2)}_{N-n}(\lambda;X_1,n_1,X_2,n_2)=
\displaystyle{\iint\limits_{\Theta}} \left( M(m_1-m_2)^+ \right.
\\
\quad + \mathbb{E}_{x}^{(2)}
R^B_{N-(n+M)}(\lambda;X_1,n_1,X_2+x,n_2+M)
\left. \right)\\
\qquad \times \lambda(m_1,m_2|X_1,n_1,X_2,n_2)dm_1 dm_2
\end{array}
\end{equation}
if $n_1+n_2<N$ where
\begin{gather*}
\mathbb{E}_{x}^{(\ell)} R(x)= \int\limits_{-\infty}^{+\infty} R(x)
f_M(x|M m_\ell) dx, \quad \ell=1,2.
\end{gather*}
Bayesian strategy prescribes currently to choose the action
corresponding to the smaller value of $R^{(1)}_{N-n}(\cdot)$,
$R^{(2)}_{N-n}(\cdot)$, the choice may be arbitrary if these
values are equal.
\par

\def \hL {\hat{L}}
\def \hR {\hat{R}}
\def \hg {\hat{g}}
\def \oX {\overline{X}}

\section{Description of the Worst-Case Prior and Corresponding Recursive Equation
}\label {WCP}

A direct usage of the main theorem of the theory of games is
virtually impossible because of the high computational complexity.
In this section, we'll specify the properties of the worst-case
prior which allow to simplify equations~\eqref{c7}-\eqref{c9}
significantly. These properties are based on the following
inequality
\begin{equation*}
R^B_N(\alpha\lambda+\tilde{\alpha}\tilde{\lambda}) \geq \alpha
R^B_N(\lambda) +\tilde{\alpha} R^B_N(\tilde{\lambda}),
\end{equation*}
if $\alpha+\tilde{\alpha}=1$; $\alpha,\tilde{\alpha}>0$, i.e.
Bayesian risk is a concave function of the prior distribution
density.

This property allows to specify the worst-case prior distribution.
Like in \cite{Koln1}, one can prove that the following
transformations $\tilde{\lambda}$ of the prior distribution
density $\lambda$ do not change the Bayesian risk, i.e.
$R^B_N(\tilde{\lambda})=R^B_N(\lambda)$:
\begin{enumerate}
\item $\tilde{\lambda}^{(1)}(m_1,m_2)= \lambda(m_2,m_1)$ (for all
$m_1$, $m_2$). This property means that expected losses do not
change if one swaps the arms of the bandit.

\item $\tilde{\lambda}^{(2)}(m_1,m_2)= \lambda(m_1+c,m_2+c)$ (for
all $m_1$, $m_2$ and any fixed $c$). This property means that
expected losses do not change if one equally shifts both
mathematical expectations.
\end{enumerate}
\par
So, if $\lambda$ is the worst-case prior distribution then
$\alpha\lambda+\tilde{\alpha}\tilde{\lambda}$ is the worst-case
prior as well. It means that the worst-case prior distribution
does not change if the above transformations are implemented. In
the sequel, it is convenient to modify parameterization. Let's put
$m_1=m+v$, $m_2=m-v$, then $\theta=(m+v,m-v)$ and
$\Theta=\{\theta:|v| \leq C\}$. Taking into account the Jacobian
$|\partial(m_1,m_2)/\partial(m,v)|=2$, a prior distribution
density is equal to $\nu(m,v)=2 \lambda(m+v,m-v)$. Then
transformations of the prior distribution densities
$\tilde{\nu}^{(1)}(m,v)=\nu(m,-v)$ and
$\tilde{\nu}^{(2)}(m,v)=\nu(m+c,v)$ (for any fixed $c$) do not
change the value of Bayesian risk. These properties allow to
specify the worst-case prior. Namely, asymptotically the
worst-case prior distribution density can be chosen the following
one:
\begin{equation}
\nu_a(m,v)=\kappa_{a}(m) \rho(v), \label{c10}
\end{equation}
where $\kappa_{a}(m)$ is the uniform density on the interval
$|m|\leq a$, $\rho(v)$ is a symmetric density (i.e.
$\rho(-v)=\rho(v)$) on the interval $|v| \leq C$  and $a \to
\infty$. This prior does not change under the first transformation
and asymptotically (as $a \to \infty$) does not change under the
second transformation.
\par
Now let's write the dynamic programming equation for calculation
the Bayesian risk with respect to~\eqref{c10}. These equations
follow from \eqref{c7}-\eqref{c9} if the prior distribution
density is formally assumed to be constant with respect to $m$ and
this gives true expressions for the posterior densities if $n_1
\geq M$, $n_2 \geq M$. At the former two steps actions should be
chosen turn-by-turn. At the time point $n=n_1+n_2$ control is
completely determined for a triple $(U,n_1,n_2)$ with $U=(X_1 n_2-
X_2 n_1)n^{-1}$.
\par
\begin{theorem} Let's put $f_D(x)=f_D(x|0)$. The strategy at the initial stage $n\le 2M$
applies actions turn-by-turn. In the sequel it can be determined
by solving the recursive Bellman-type equation:
\begin{gather}
\label{c11}
R_M(U,n_1,n_2)=\min_{\ell=1,\,2}R^{(\ell)}_M(U,n_1,n_2),
\end{gather}
where $R^{(1)}_M(U,n_1,n_2)=R^{(2)}_M(U,n_1,n_2)=0$ if $n_1+n_2=N$
and
\begin{gather}
\label{c12} R^{(1)}_M(U,n_1,n_2)=M g^{(1)}(U,n_1,n_2)\\+
\displaystyle{\int\limits_{-\infty}^{\infty}} R_M(U-x,n_1+M,n_2)
f_{M n_2^2 n^{-1} (n+M)^{-1}}(x) dx,\nonumber
\\R^{(2)}_M(U,n_1,n_2)=M g^{(2)}(U,n_1,n_2)\\+ \label{c13}
\displaystyle{\int\limits_{-\infty }^{\infty}} R_M(U-x,n_1,n_2+M)
f_{M n_1^2 n^{-1} (n+M)^{-1}}(x) dx \nonumber
\end{gather}
if $n_1+n_2<N$. Here
\begin{gather}\label{c14}
\begin{array}{c}
g^{(\ell)}(U,n_1,n_2)\\= \displaystyle{\int\limits_0^C} 2v
\exp\left((-1)^{\ell}2Uv -2v^2 n_1 n_2 n^{-1} \right)\rho(v) dv,
\end{array}
\end{gather}
$\ell=1,2$. If $n> 2M$ then the $\ell$-th action is currently
optimal iff $R^{(\ell)}_M(U,n_1,n_2)$ has smaller value
($\ell=1,2$). Corresponding Bayesian risk~\eqref{c4} is calculated
as follows
\begin{gather}
\label{c15}
\lim\limits_{a \to \infty} R^B_N(\nu_a(m,v))= R^B_N(\rho(v))\\
=4 M \displaystyle{\int\limits_0^C} v \rho(v) dv+
\displaystyle{\int\limits_{-\infty}^{\infty}} R_M(U,M,M) f_{0.5 M}
(U) dU. \nonumber
\end{gather}
\end{theorem}
\par
Proof of theorem is presented in Appendix A.

\par
\section{Invariant Recursive Equation and Passage to the Limit } \label {IRE}
\par
Let's introduce the following change of variables $\eps=M N^{-1}$,
$t_1=n_1N^{-1}$, $t_2=n_2N^{-1}$, $t=nN^{-1}$, $u=U N^{-1/2}$,
$w=vN^{1/2}$, $c=CN^{1/2}$, $\varrho(w) = N^{1/2}\rho(v)$,
$r_\eps(u,t_1,t_2)=N^{-1/2}R_M(U,n_1,n_2)$,\\
$r^{(\ell)}_\eps(u,t_1,t_2)=N^{-1/2}R^{(\ell)}_M(U,n_1,n_2)$. Now
we consider the set of close expectations
\begin{gather*}
\Theta_N=\{w:w \leq c\}=\{\theta:|m_1-m_2|\leq 2cN^{-1/2}\}.
\end{gather*}
Recall that according to~\eqref{c3} the maximal expected losses in
the two-armed bandit problem have the order $N^{1/2}$ and are
attained just for close expectations with $c>0$ large enough. On
the contrary, the maximal expected losses for distant expectations
$|m_1-m_2|\ge \delta>0$ have the order $\log(N)$. This estimate
follows from the results of~\cite{Lai}.

For close expectations the following theorem holds.
\par
\begin{theorem} The strategy at the initial stage $t\le 2\eps$ $(n \le
2\eps N)$ applies actions turn-by-turn. Then it can be determined
by solving the following recursive Bellman-type equation:
\begin{gather}
\label{c20}
r_\eps(u,t_1,t_2)=\min_{\ell=1,\,2}r^{(\ell)}_\eps(u,t_1,t_2),
\end{gather}
where $r^{(1)}_\eps(u,t_1,t_2)=r^{(2)}_\eps(u,t_1,t_2)=0$ if
$t_1+t_2=1$ and
\begin{gather}
\label{c21} r^{(1)}_\eps(u,t_1,t_2)=\eps g^{(1)}(u,t_1,t_2)\\
+ \displaystyle{\int\limits_{-\infty}^{\infty}}
r_\eps(u-x,t_1+\eps,t_2) f_{\eps t_2^2 t^{-1} (t+\eps)^{-1}}(x)
dx,\nonumber
\\r^{(2)}_\eps(u,t_1,t_2)=\eps g^{(2)}(u,t_1,t_2) \label{c22}
\\+
\displaystyle{\int\limits_{-\infty }^{\infty}}
r_\eps(u-x,t_1,t_2+\eps) f_{\eps t_1^2 t^{-1} (t+\eps)^{-1}}(x) dx
\nonumber
\end{gather}
if $t_1+t_2<1$. Here
\begin{gather}\label{c23}
g^{(\ell)}(u,t_1,t_2)\\= \displaystyle{\int\limits_0^c} 2w
\exp\left((-1)^{\ell}2uw -2w^2 t_1 t_2 t^{-1} \right)\varrho(w)
dw,\nonumber
\end{gather}
$\ell=1,2$. If $t> 2\eps$ $(n> 2\eps N)$ then the $\ell$-th action
is currently optimal iff $r^{(\ell)}_\varepsilon(u,t_1,t_2)$ has
smaller value ($\ell=1,2$). Bayesian risk corresponding to the
worst-case prior distribution is calculated according to the
formula
\begin{gather}
\label{c24} N^{-1/2}R^B_N(\rho(v))\\
=4 \eps \displaystyle{\int\limits_0^c} w \varrho(w) dw +
\int\limits_{-\infty}^{\infty} r_\eps(u,\eps,\eps) f_{0.5 \eps}
(u) du \nonumber
\end{gather}
\end{theorem}
{\bf Proof.} Is done by implementation of described above change
of variables.
\par
Let's denote by $r_\eps(\varrho;u,t_1,t_2)$ the Bayesian risk as
dependent on a prior distribution $\varrho(w)$. Obviously,
$r_\eps(\varrho;u,t_1,t_2)$ is a decreasing function of $\eps$ for
any fixed $u$, $t_1$, $t_2$ because diminishing of $\eps$ implies
that actions may be changed more often. The following theorem is
given without proof.
\begin{theorem}\label{th3}
For all  $u$, $t_1$, $t_2$, for which the solution to equation
\eqref{c20}-\eqref{c22} is well defined, there exist limits
$r(\varrho;u,t_1,t_2)=\lim\limits_{\eps \to 0}
r_\eps(\varrho;u,t_1,t_2)$ which can be extended by continuity to
all $u,\, t_1,\, t_2\, (t_1>0,\, t_2>0,\, t_1+t_2<1)$. These
limits are uniformly bounded and satisfy Lipschitz conditions in
$u$. The minimax risk on the set of close expectations
$\Theta_N=\{|m_1-m_2|\leq 2cN^{-1/2}\}$ satisfies the equality
\begin{equation}\label{c25}
\lim_{N \to \infty} N^{-1/2} R^M_N(\Theta_N) = \sup_\varrho
r(\varrho;0,0,0),
\end{equation}
where $r(\varrho;0,0,0)=\lim\limits_{\eps \to
0}r(\varrho;0,\eps,\eps)$.
\end{theorem}
\par
\par
Let's present the limiting description of $r(u,t_1,t_2)$ by the
second order partial differential equation. Assume that
$r_\eps(u,t_1,t_2)$ has continuous partial derivatives of proper
orders. We present $r_\eps(u-x,t_1+\eps,t_2)$ as Taylor series:
\begin{gather}
\label{c26}
\begin{array}{c}
r_\eps(u-x,\cdot)=r_\eps(u,\cdot)-x
\times\displaystyle{\frac{\partial r_\eps(u,\cdot)}{\partial
u}}\\+ \displaystyle{\frac{x^2}{2}} \times \frac{\partial^2
r_\eps(u,\cdot)}{\partial u^2}+o(x^2).
\end{array}
\end{gather}
Noting that
\begin{gather*}
\begin{array}{c}
\displaystyle{\int_{-\infty}^\infty} f_\eps (x) dx=1,
\displaystyle{\int_{-\infty}^\infty} x f_\eps(x) dx= 0,
\displaystyle{\int_{-\infty}^\infty} x^2 f_\eps(x) dx= \eps,
\end{array}
\end{gather*}
and substituting \eqref{c26} into \eqref{c21} one obtains
\begin{gather}\label{c27}
\begin{array}{l}
r^{(1)}_\eps(u,t_1,t_2)=\eps g^{(1)}(u,t_1,t_2)\\+
\displaystyle{\int\limits_{-\infty}^\infty}
r_\eps(u-x,t_1+\eps,t_2)
f_{\eps t_2^2 t^{-1} (t+\eps)^{-1}}(x) dx\\
=\eps g^{(1)}(u,t_1,t_2)+
r_\eps(u,t_1+\eps,t_2)\\+\displaystyle{\frac{\eps t_2^2}{2 t
(t+\eps)}}\times \frac{\partial^2 r_\eps(u,t_1+\eps,t_2)}{\partial
u^2}+o(\eps),
\end{array}
\end{gather}
Similarly,
\begin{gather}
\begin{array}{l}\label{c28}
r^{(2)}_\eps(u,t_1,t_2)=\eps g^{(2)}(u,t_1,t_2)+r_\eps(u,t_1,t_2+\eps)\\
\\+\displaystyle{\frac{\eps t_1^2}{2 t (t+\eps)}} \times
\frac{\partial^2 r_\eps(u,t_1,t_2+\eps)}{\partial u^2}+o(\eps).
\end{array}
\end{gather}
Recall now that equations \eqref{c27}-\eqref{c28} must be
complemented by equation \eqref{c20} which can be written as
\begin{gather}\label{c29}
\min_{\ell=1,\,2}
(r^{(\ell)}_\eps(u,t_1,t_2)-r_\eps(u,t_1,t_2))=0.
\end{gather}
From~\eqref{c27}-\eqref{c29} one obtains (as $\eps \downarrow 0$)
the partial differential equation:
\begin{gather} \label{c30}
\min_{\ell=1,\,2}\left(\frac{\partial r}{\partial t_\ell} +
\frac{t_{\overline{\ell}}^2}{2 t^2}\times \frac{\partial^2
r}{\partial u^2} +g^{(\ell)}(u,t_1,t_2)\right)=0
\end{gather}
with $\overline{\ell}=3-\ell$. Initial and boundary conditions
take the form
\begin{gather}\label{c31}
\begin{array}{l}
r(u,t_1,t_2)\|_{t_1+t_2=1}=0, \\
r(\infty,t_1,t_2)= r(-\infty,t_1,t_2)=0.
\end{array}
\end{gather}
The optimal strategy prescribes to chose the $\ell$-th action if
the the $\ell$-th member in the left-hand side of~\eqref{c30} has
minimal value.
\par
\section{Numerical Results } \label{CIS}
\par
\begin{figure}
\begin{center}
\includegraphics[width=8.4cm]{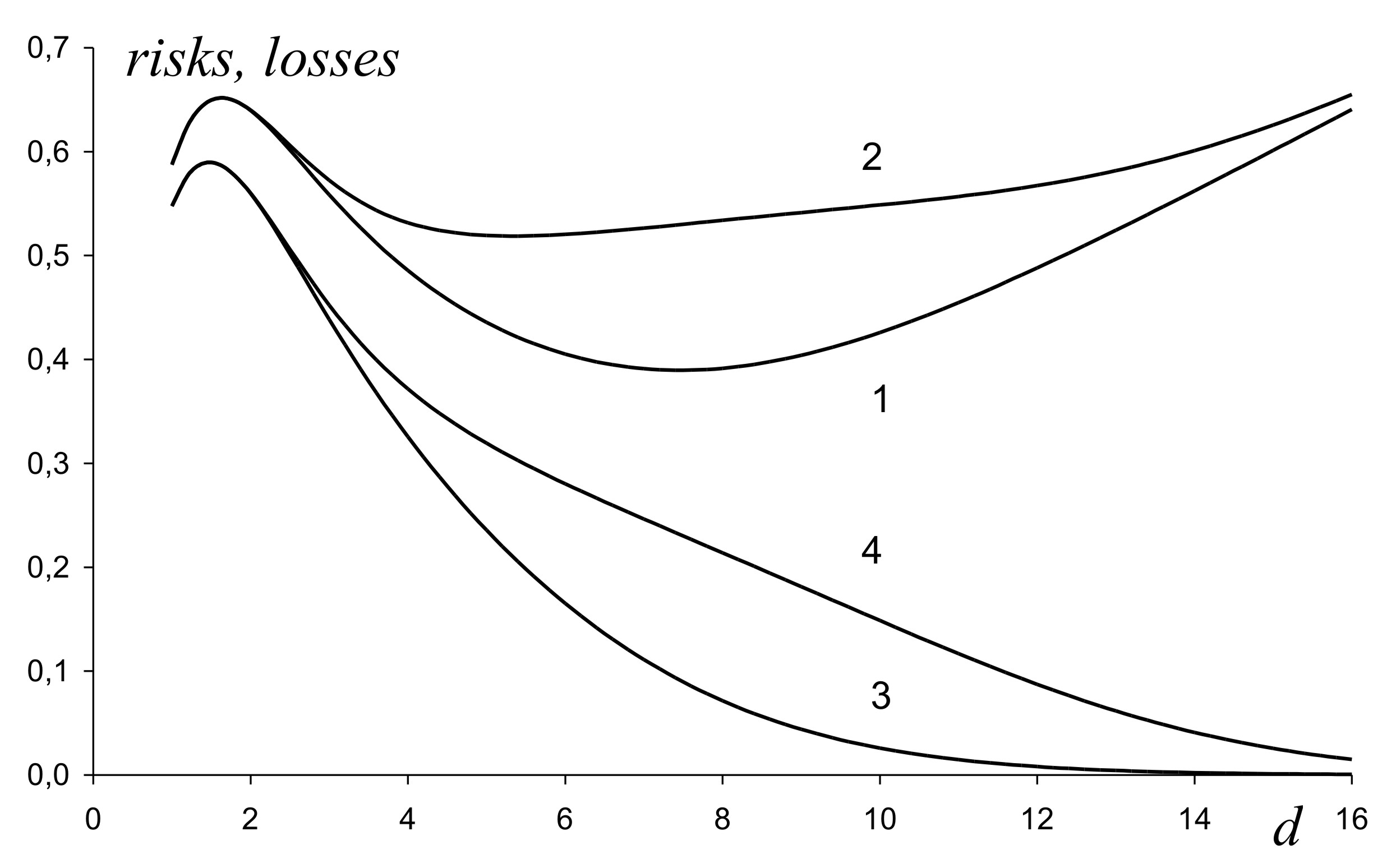} 
\caption{Risks, expected losses} \label{fig1}
\end{center}
\end{figure}

Bayesian risks were calculated by~\eqref{c20}-\eqref{c24} with
$\eps=0.02$. It was assumed that the worst-case prior $\varrho(w)$
is a degenerate one and concentrated at two points $w=\pm d$ with
equal probabilities 0.5. The risks are presented by line 1 on
figure~\ref{fig1} as a function of $d$. The worst-case prior
corresponds to its maximum. The maximum is approximately equal to
$0.65$ at $d \approx 1.6$.
\par
Expected losses corresponding to determined strategy
$\sigma_\ell(u,t_1,t_2)=\Pr(y_n=\ell|u,t_1,t_2)$ were sought for
by solving recursive equation
\begin{gather*}
l_\eps(u,t_1,t_2)=\sigma_1(u,t_1,t_2)l^{(1)}_\eps(u,t_1,t_2)\\+
\sigma_2(u,t_1,t_2)l^{(2)}_\eps(u,t_1,t_2),
\end{gather*}
where
\begin{gather*}
l^{(1)}_\eps(u,t_1,t_2)=l^{(2)}_\eps(u,t_1,t_2)=0
\end{gather*}
if $t_1+t_2=1$ and then
\begin{gather*}
l^{(1)}_\eps(u,t_1,t_2)=\eps g^{(1)}(u,t_1,t_2)\\
+ \displaystyle{\int\limits_{-\infty}^{\infty}}
l_\eps(u-x,t_1+\eps,t_2) f_{\eps t_2^2 t^{-1} (t+\eps)^{-1}}(x)
dx,\nonumber
\\l^{(2)}_\eps(u,t_1,t_2)=\eps g^{(2)}(u,t_1,t_2)
\\+
\displaystyle{\int\limits_{-\infty }^{\infty}}
l_\eps(u-x,t_1,t_2+\eps) f_{\eps t_1^2 t^{-1} (t+\eps)^{-1}}(x) dx
\nonumber
\end{gather*}
if $t_1+t_2<1$. Then
\begin{gather*}
N^{-1/2}L_N(\sigma,\rho(v))\\
=4 \eps
\displaystyle{\int\limits_0^c} w \varrho(w) dw +
\int\limits_{-\infty}^{\infty} l_\eps(u,\eps,\eps) f_{0.5 \eps}
(u) du.
\end{gather*}
The losses are presented by line 2 on figure~\ref{fig1}. One can
see that its maximal value does not exceed the value 0.65 and this
confirms the assumption concerning the worst-case prior.
Nevertheless, one can see that expected losses become larger than
0.65 if $d>16$. This is caused by the initial stage of control
where both actions are equally applied. On figure~\ref{fig1} lines
3 and 4 present risks and expected losses without those ones at
the initial stage. These functions do not grow with growing $d$.
Therefore, to reduce expected losses at large $d$ one should
reduce initial stage of control.
\par
To obtain the limiting value of the minimax risk~\eqref{c25}
calculations of $r(\varrho;u,t_1,t_2)$, as a function of $d$, were
implemented according to~\eqref{c20}, \eqref{c27}, \eqref{c28},
\eqref{c31} with $\eps=0.001$ for $|u| \le 2.3$. Partial
derivatives were replaced by partial differences with $\Delta
u=0.023$, $\Delta t=2000^{-1}$. It was assumed that $\varrho(w)$
is a degenerate distribution density concentrated at two points
$w=\pm d$. For $0.5\le d \le 2.5$ maximum of $2d\eps +
r(\varrho;0,\eps,\eps)$ was approximately equal to 0.637 at $d
\approx 1.57$. Hence, the minimax risk corresponding to batch
processing in 50 stages is approximately 2\% larger than the
limiting value.
\par
Monte-Carlo simulations were implemented for batch processing of
$T=5000$ items of data by packets of $M=100$ data items, i.e. in
50 stages. The normalized expected losses $(DT)^{-1/2} L_T(\sigma,
\theta)$ with $\theta=(p+d(D/T)^{1/2},p-d(D/T)^{1/2})$, $p=0.5$,
$D=0.25$ were calculated as a function of $d$. This function is
just the same as the line 2 on figure 1 and that is why it is not
specially presented there.
\par
\section{Conclusion } \label{Con}
\par
The minimax approach to the two-armed bandit problem based on the
main theorem of the theory of games is proposed. Incomes of the
two-armed bandit are assumed to have Gaussian distributions and
this implies the possibility of their batch processing. The
approach allows to determine numerically minimax strategy and
minimax risk for any finite control horizon by solving
Bellman-type recursive equation. However, the results have an
asymptotic nature because they should be applied to batch
processing a large amount of data by packets in a moderate number
of stages. At the initial stage of control, there may be large
expected losses because at initial stage actions are chosen
turn-by-turn. To reduce losses at the initial stage one should
take initial packets of data having smaller sizes.


\begin{thebibliography}{xx}  

\bibitem [Bather (1983)] {Bather}
Bather, J. A. (1983).
\newblock The minimax risk for the two-armed bandit
problem.
\newblock \emph{ Lecture Notes in Statistics},
\newblock volume 20, 1--11. Springer-Verlag, New York.

\bibitem[Berry and Fristedt (1985)] {BF}
Berry, D. A., and Fristedt, B. (1985).
\newblock \emph{  Bandit Problems: Sequential
Allocation of Experiments}.
\newblock Chapman and Hall, London, New~York.

\bibitem [Cheng (2003)] {Cheng}
Cheng, T., Su, Y., Berry, D.A. (2003)
\newblock Choosing sample size for a clinical trial using decision
analysis.
\newblock \emph{Biometrika.} 90, 923--936.

\bibitem [Fabius and van Zwet (1970)] {FZ}
Fabius, J., and van~Zwet, W. R. (1970).
\newblock Some remarks on the two-armed bandit.
\newblock \emph{Ann. Math. Statist.}, 41, 1906--1916.

\bibitem [Gasnikov et al (2015)] {GNS}
Gasnikov, A. V., Nesterov, Yu. E., and Spokoiny, V. G. (2015).
\newblock On the efficiency of a randomized mirror descent algorithm in
online optimization problems.
\newblock \emph{ Computational Mathematics and Mathematical Physics}, 55:\penalty0 4,
580--596.

\bibitem [Juditsky et al (2008)] {UN}
Juditsky, A., Nazin, A. V., Tsybakov, A. B., and Vayatis, N.
(2008).
\newblock Gap-free bounds for stochastic multi-armed bandit.
\newblock \emph{  Proc. 17th World Congress IFAC}. Seoul, Korea, July 6--11),
11560--11563.

\bibitem[Kolnogorov (2010)] {Koln1}
Kolnogorov, A. V. (2010).
\newblock
Determination of the minimax risk for the normal two-armed bandit.
\newblock
In \emph{Proceedings of the IFAC Workshop ``Adaptation and
Learning in Control and Signal Processing ALCOSP 2010''}, Antalya,
Turkey, August 26--28, 2010. DOI 10.3182/20100826-3-TR-4015.00044.
http://www.ifac-papersonline.net.

\bibitem [Kolnogorov (2012)]{Koln2}
Kolnogorov, A. V. (2012).
\newblock Parallel design of robust control in the stochastic
environment (the two-armed bandit problem).
\newblock \emph{   Automation and Remote Control}, 73:\penalty0 4, 689--701.

\bibitem [Kolnogorov (2015)] {Koln3}
Kolnogorov, A. V. (2015).
\newblock On a limiting description of robust
parallel control in a random environment.
\newblock \emph{  Automation and Remote Control}, 76:\penalty0 7, 1229--1241.

\bibitem [Lai et al (1980)] {Lai}
Lai, T. L., Levin, B., Robbins, H., and Siegmund, D. (1980).
\newblock Sequential medical trials (stopping rules/asymptotic
optimality.
\newblock \emph{ Proc. Nati. Acad. Sci. USA.} 77:\penalty0 6, 3135--3138.

\bibitem [Lugosi and Cesa-Bianchi (2006)] {Lugosi}
Lugosi, G., and Cesa-Bianchi, N. (2006).
\newblock \emph{ Prediction, Learning
and Games}.
\newblock  Cambridge University Press, New York.

\bibitem [Nazin and Poznyak (1986)] {NP}
Nazin, A. V., and Poznyak, A. S. (1986).
\newblock \emph{Adaptive Choice of Alternatives}.
\newblock Nauka, Moscow. (In Russian)

\bibitem [Presman and Sonin (1990)] {PS}
Presman, E. L., and Sonin, I. M. (1990).
\newblock \emph{Sequential Control with Incomplete Information}.
\newblock  Academic Press, New York.

\bibitem [Robbins (1952)] {Robbins}
Robbins, H. (1952).
\newblock Some aspects of the sequential design of experiments.
\newblock \emph{ Bulletin AMS.}, 58:\penalty0 5, 527--535.

\bibitem [Sragovich (2006)]{Sragovich}
Sragovich, V. G. (2006).
\newblock \emph{Mathematical Theory of Adaptive
Control}.
\newblock World Scientific. Interdisciplinary Mathematical Sciences,  New
Jersey, London,  volume 4.

\bibitem [Tsetlin (1973)] {Tsetlin}
Tsetlin, M. L. (1973).
\newblock \emph{Automation Theory and Modeling of Biological Systems}.
\newblock Academic Press, New York.

\bibitem [Varshavsky (1973)] {Varsh}
Varshavsky, V. I. (1973).
\newblock \emph{Collective Behavior of
Automata}.
\newblock  Nauka, Moscow. (In Russian)

\bibitem [Vogel (1960)] {Vogel}
Vogel, W. (1960).
\newblock An asymptotic minimax theorem for the
two-armed bandit problem.
\newblock \emph{Ann. Math. Stat.}, 31, 444--451.








\end{thebibliography}

\appendix
\section{Proof of Theorem 2}    

{\bf Proof.} Let's put
\begin{gather*}
\hR(X_1,n_1,X_2,n_2)\\=R^B_{N-n}(X_1,n_1,X_2,n_2)p(X_1,n_1,X_2,n_2)
\end{gather*}
with $p(X_1,n_1,X_2,n_2)$ defined in~\eqref{c6}.  Denote by\\
$\hR(Z,n_1,n_2)=\hR(X_1,n_1,X_2,n_2)$ with $Z=X_1 n_2- X_2 n_1$.
Let's check that if the prior is given by~\eqref{c10}
then~\eqref{c8} takes the form
\begin{gather}\label{c16}
\hR^{(1)}(Z,n_1,n_2)=\int\limits_0^C 2 M v \hg(Z,n_1,n_2,v)\rho(v)
dv
\\
+n_2^{-1}\int\limits_{-\infty}^{+\infty}
\hR(Z+z,n_1+M,n_2)h_{M}\left(\frac{MZ-n_1z}{n_2},n_1\right)
dz,\nonumber
\end{gather}
with
\begin{equation}
\begin{array}{ll}\label{c17}
\hg(Z,n_1,n_2)
=\displaystyle{\frac{1}{(2\pi n_1 n_2(n_1+n_2))^{1/2}}}  \\
\qquad \times \exp\left(-\displaystyle{\frac{(Z + 2vn_1 n_2)^2}{2
n_1 n_2 (n_1+n_2)}}\right),
\end{array}
\end{equation}
and
\begin{equation}\label{c18}
h_{M}(z,n)= \left(\frac{n+M}{2\pi M n}\right)^{1/2} \times \exp
\left( -\frac{z^2}{2Mn(n+M)}\right).
\end{equation}
\par
Really, if the prior is taken from~\eqref{c10} then \eqref{c8}
takes the form
\begin{gather}\nonumber
\hR^{(1)}(X_1,n_1,X_2,n_2)\\
\label{c19}=\displaystyle{\int\limits_0^C} 2 M v
\hg(X_1,n_1,X_2,n_2,v)\rho(v) dv
\\
+\displaystyle{\int\limits_{-\infty}^{+\infty}}
\hR(X_1+x,n_1+M,X_2,n_2)h_{M}\left(MX_1-n_1x,n_1\right) dz.
\nonumber
\end{gather}
Here
\begin{equation*}
\begin{array}{ll}
\hg(X_1,n_1,X_2,n_2)\\
=\displaystyle{\int_{-\infty}^{\infty}}f_{n_1}(X_1|n_1(m+v)) f_{n_2}(X_2|n_2(m-v))dm\\
=\displaystyle{\frac{1}{(2\pi n_1 n_2(n_1+n_2))^{1/2}}}  \times
\exp\left(-\displaystyle{\frac{(\oX_1 -\oX_2 + 2v)^2}{2
(n_1^{-1}+n_2^{-1})}}\right),
\end{array}
\end{equation*}
and
\begin{equation*}
\begin{array}{l}
h_{M}(MX_1-n_1 x,n_1)= \\
=\Big(  \displaystyle{\iint\limits_{\Theta}}  f_M(x|M m_1)
f_{n_1}(X_1|n_1 m_1) \times \\
\times f_{n_2}(X_2|n_2 m_2) \lambda(m_1,m_2) dm_1 dm_2 \Big) \Big/ \\
\Big/ \Big( \displaystyle{\iint\limits_{\Theta}}
f_{n_1+M}(X_1+x|(n_1+M)
m_1) \times \\
f_{n_2} (X_2|n_2)
m_2) \lambda(m_1,m_2) dm_1 dm_2 \Big) = \\
=\displaystyle{\frac {f_M(x|Mm_\ell) f_{n_1}(X_1|n_1 m_1)} {
f_{n_1+M}(X_1+x|(n_1+M)m_1)}} =\\
=\left(\displaystyle{\frac{n_1+M}{2\pi n_1 M}}\right)^{1/2}\exp
\left( -\displaystyle{\frac{(MX_1-n_1 x)^2}{2 n_1 M
(n_1+M)}}\right).
\end{array}
\end{equation*}
So, these expressions correspond to~\eqref{c17}-\eqref{c18}. Note
that at $n_1+M$, $n_2$ the value $Z$ is recalculated by expression
$Z \leftarrow (X_1 + x)n_2 - X_2(n_1 +M) = Z +z$ with $z = xn_2
-MX_2$. Noting that $MX_1 - n_1x = n_2^{-1}(ZM-n_1z)$ and changing
the integration variable in \eqref{c19} from $x$ to $z$ one
obtains \eqref{c16}.
\par
Now let's put $\hR(Z,n_1,n_2)=f_{n_1 n_2 n}(Z) R(U,n_1,n_2)$. The
first equation \eqref{c12} may be obtained from \eqref{c16} and
equality
\begin{gather*}
\frac{n_2^{-1}h_M\left(n_2^{-1}(Z(n_1 + M) -n_1
y),n_1\right)}{f_{n_1 n_2 n} (Z)}\\  = \frac{ f_{M n_2^2
(n+M)^{-1}n^{-1}}(y(n+M)^{-1}-U)}{(n+M)\times f_{(n_1+M) n_2
(n+M)} (y)}.
\end{gather*}
The second equation \eqref{c12} is similarly checked. Obviously,
Bayesian risk \eqref{c4} is calculated according to the formula
\begin{equation*}
R^B_N(\rho(v))= 4M\displaystyle{\int\limits_0^C} v \rho(v) dv+
\int\limits_{-\infty}^{\infty} \hR(z,M,M) dz
\end{equation*}
and hence by \eqref{c15}.

\end{document}